\documentclass[oneside,english]{amsart}
\usepackage[T1]{fontenc}
\usepackage[latin9]{inputenc}
\usepackage{amsthm}
\usepackage{amstext}
\usepackage{amssymb}
\usepackage{esint}
\usepackage{babel}

\makeatletter

\newcommand{\noun}[1]{\textsc{#1}}

\numberwithin{equation}{section}
\numberwithin{figure}{section}
\theoremstyle{plain}
\newtheorem{thm}{Theorem}[section]
  \theoremstyle{plain}
  \newtheorem{prop}[thm]{Proposition}
  \theoremstyle{definition}
  \newtheorem{defn}[thm]{Definition}
  \theoremstyle{plain}
  \newtheorem{lem}[thm]{Lemma}
  \theoremstyle{plain}
  \newtheorem{cor}[thm]{Corollary}
 \theoremstyle{definition}
  \newtheorem{example}[thm]{Example}
  \theoremstyle{remark}
  \newtheorem{rem}[thm]{Remark}

\makeatother

\begin{document}

\title{$L^{2}$-estimates for the $d$-operator acting on super forms}

\author{Aron Lagerberg}
\begin{abstract}
In the setting of super forms developed in \cite{Lagerberg}, we introduce
the notion of $\mathbb{R}-$K{\"a}hler metrics on $\mathbb{R}^{n}$. We
consider existence theorems and $L^{2}-$estimates for the equation
$d\alpha=\beta$, where $\alpha$ and $\beta$ are super forms, in
the spirit of H{\"o}rmander's $L^{2}-$estimates for the $\bar{\partial}-$equation
on a complex K{\"a}hler manifold. 
\end{abstract}

\address{A. Lagerberg: Department of Mathematics, Chalmers University of Technology
and the University of G{\"o}teborg, 412 96, G{\"o}TEBORG, SWEDEN.}

\email{aronl@chalmers.se}

\maketitle
\tableofcontents{}

\section{Introduction}

This article is concerned with introducing the notion of an $\mathbb{R}-$K{\"a}hler
metric on the Euclidean space, $\mathbb{R}^{n}$. Let us explain the
meaning of this statement: on a complex manifold, a hermitian metric
induces a $(1,1)-$form $\omega$, and the manifold is \emph{K{\"a}hler}
if $d\omega=0$. In \cite{Lagerberg} the formalism of super forms
on $\mathbb{R}^{n}$ was considered, which enables us to define $(p,q)-$forms
on $\mathbb{R}^{n}$. In particular, a smooth metric $g$ on $\mathbb{R}^{n}$
can be represented by a smooth, positive $(1,1)-$form $\omega$,
and in analogy with the complex setting, we define the metric $g$
to be \emph{$\mathbb{R}-$K{\"a}hler} if $d\omega=0$. In this article,
our main concern is for the $d-$equation for $(p,q)$-forms on $\mathbb{R}^{n}$
endowed with a K{\"a}hler metric; by this we mean that given a $(p,q)$-form
$\beta$, we wish to find a $(p-1,q)$-form $\alpha$ solving the
equation \[
d\alpha=\beta.\]
Under certain hypothesis on $\beta$, we shall prove existence theorems
for this equation using arguments from the technique of $L^{2}-$estimates
due to H{\"o}rmander for the $\overline{\partial}-$equation on a complex
K{\"a}hler manifold (c.f. \cite{Hormander}). This will also give us an
$L^{2}-$estimate on the solution $\alpha$ in terms of $\beta$ on
a given $L^{2}$-space (depending on the K{\"a}hler metric), to be introduced
later in this article. As a particular case, we are able to solve
the $d-$equation for ordinary $p$-forms on $\mathbb{R}^{n}$ together
with an $L^{2}-$estimate on the solution in terms of the given data.
The key point in applying the arguments of H{\"o}rmander is to establish
a Kodaira-Bochner-Nakano-type identity (c.f \cite{Nakano}) for natural
Laplace-operators arising in our setting. We also take the opportunity
to introduce, in analogy with the complex case, the theory of primitive
super forms. Our hope is that the results developed in this article
can be used to establish results in convex analysis. For instance,
there are many articles concerned with convex inequalities that utilizes
$L^{2}$-theory (see for instance \cite{Cordero},\cite{Berndtsson2}),
and we hope that our approach in this article will give a fruitful
addition to the theory already developed.

\medskip{}

\noun{\textbf{{Acknowledgements:}}} I would like to thank my advisor
Bo Berndtsson for inspiration and support.

\section{Preliminaries\label{sec:Preliminaries}}

In this article, we will consider differential forms in $\mathbb{R}^{n}\times\mathbb{R}^{n}=\{(x_{1},...,x_{n},\xi_{1},...,\xi_{n})\}$
with coefficients depending only on the variables $(x_{1},...,x_{n})$.
Such forms, which we shall call super forms, were considered in the
article \cite{Lagerberg}. We say that $\alpha$ is a $(p,q)-$form
if \[
\alpha=\sum_{|I|=p,|J|=q}\alpha_{IJ}(x)dx_{I}\wedge d\xi_{J},\]
where we use multi-index notation, and a $k$-form is a $(p,q)$-form
with $p+q=k$. The set of $(p,q)-$forms whose coefficients are smooth
will be denoted by $\mathcal{E}^{p,q}$. A smooth $(1,1)$-form $\omega$
is said to be \emph{positive} if the coefficient-matrix $(\omega_{ij}(x))_{i,j}$
is positive definite, for each $x$. Let us define the operator \[
d^{\#}:\mathcal{E}^{p,q}\rightarrow\mathcal{E}^{p,q+1}\]
by letting \[
d^{\#}(\sum_{|I|=p,|J|=q}\alpha_{IJ}(x)dx_{I}\wedge d\xi_{J})=\sum_{i=1}^{n}\sum_{|I|=p,|J|=q}\frac{\partial}{\partial x_{i}}\alpha_{IJ}(x)d\xi_{i}\wedge dx_{I}\wedge d\xi_{J}.\]
The operator $d:\mathcal{E}^{p,q}\rightarrow\mathcal{E}^{p+1,q}$
is defined as usual, and a form $\alpha$ is called \emph{closed}
if $d\alpha=0$. We also define the linear map \[
J:\{(p,q)-\textmd{forms}\}\rightarrow\{(q,p)-\textmd{forms}\}\]
by letting \[
J(\sum_{|I|=p,|J|=q}\alpha_{IJ}(x)dx_{I}\wedge d\xi_{J})=\sum_{|I|=p,|J|=q}\alpha_{IJ}(x)d\xi_{I}\wedge dx_{J}.\]
Observe that this makes for $J^{2}=Id$. The operator $d^{\#}$ can
be written in terms of $J$ as \[
d^{\#}=J\circ d\circ J.\]
 We have the following result (cf. \cite{Lagerberg}):
\begin{prop}
A closed, smooth $(1,1)-$form $\omega$ is positive if and only if
there exists a convex, smooth function $f$ such that \[
\omega=dd^{\#}f.\]
 
\end{prop}
Now fix a smooth, positive, and closed $(1,1)$- form $\omega$. We
shall use the notation \[
\omega_{q}=\omega^{q}/q!.\]
 Such a form $\omega$ induces a metric on $\mathbb{R}^{n}$ in a
natural way: if $v=(v_{1},...,v_{n}),w=(w_{1},...,w_{n})\in\mathbb{R}^{n}$,
then for every $x\in\mathbb{R}^{n}$, we define\[
(v,w)_{x}=\sum_{i,j=1}^{n}v_{i}w_{j}\omega_{ij}(x),\]
where the functions $\omega_{ij}$ are defined by $\omega=\sum_{i,j=1}^{n}\omega_{ij}dx_{i}\wedge d\xi_{j}$.
We obtain an induced metric on the space of $(1,0)-$ and $(0,1)-$forms:
if $\alpha=\sum\alpha_{i}dx_{i}$ then $(\alpha,\alpha)_{x}=\sum\omega^{ij}(x)\alpha_{i}\alpha_{j}$
where $(\omega^{ij})$ denotes the inverse of the matrix $(\omega_{ij})$,
and analogously for $(0,1)$-forms. Using this metric, we would like
to define the norm of a $(p,q)$-form, at a point. Let us fix an orthonormal
(with respect to $\omega$) coordinate system $(dx_{1},...,dx_{n})$
for the space of $(1,0)$- forms. If $\alpha=\sum\alpha_{IJ}dx_{I}\wedge d\xi_{J}$,
we define \begin{equation}
|\alpha|^{2}=\sum|\alpha_{IJ}|^{2}.\label{eq:definitionofnormofform}\end{equation}
If $\alpha=\sum_{|I|=p}\alpha_{I}dx_{I}$ and $\beta=\sum_{|J|=q}\beta_{J}d\xi_{J}$,
then \[
|\alpha\wedge\beta|^{2}=\sum_{I,J}\alpha_{I}^{2}\beta_{J}^{2}=|\alpha|^{2}|\beta|^{2}.\]
If we polarize this formula we obtain\begin{equation}
(\alpha\wedge\beta,\alpha'\wedge\beta')=(\alpha,\alpha')(\beta,\beta'),\label{eq:polarizednorm}\end{equation}
with $(p,0)-$forms $\alpha,\alpha'$ and $(0,q)-$forms $\beta$,$\beta'$,
and where $(\cdot,\cdot)$ denotes the inner product associated with
the norm $|\cdot|$. Let us show that the definition \eqref{eq:definitionofnormofform}
is independent of the choice of orthonormal coordinate system: We
begin with the case of $(p,0)$- forms: Let $\alpha=\sum_{|I|=p}\alpha_{I}dx_{I}.$
A simple calculation shows that, \[
|\alpha|^{2}\omega_{n}=c_{p}\alpha\wedge J(\alpha)\wedge\omega_{n-p},\]
where $c_{p}=(-1)^{p(p-1)/2}$, and this expression does not depend
on the basis chosen. The number $c_{p}$ is chosen such that $dx_{i_{1}}\wedge...\wedge dx_{i_{p}}\wedge d\xi_{i_{1}}\wedge...\wedge d\xi_{i_{p}}=c_{p}\cdot dx_{i_{1}}\wedge d\xi_{i_{1}}\wedge...\wedge dx_{i_{p}}\wedge d\xi_{i_{p}}$.
The same calculations hold for $(0,q)$-forms. Thus, at least for
$(p,0)$- and $(0,q)-$forms, formula \eqref{eq:definitionofnormofform}
does not depend on which orthonormal coordinates we choose. Now, let
$(dy_{1},...,dy_{n})$ be another orthonormal basis, and let $d\zeta_{i}=J(dy_{i})$.
If $\alpha=\sum\alpha_{IJ}dy_{I}\wedge d\zeta_{J}$, then by \eqref{eq:polarizednorm}
we get that \[
(\alpha,\alpha)=\sum\alpha_{IJ}\alpha_{KL}(dy_{I},dy_{K})(d\zeta_{J},d\zeta_{L}).\]
But by the above, we know that $(\cdot,\cdot)$ does not depend on
which orthonormal basis we work with, when applied to $(p,0)-$ or
$(0,q)-$forms. Thus $(dy_{I},dy_{K})$ and $(d\zeta_{J},d\zeta_{L})$
is non zero, and equal to one, if and only if $I=K$ and $J=L$. Thus
the definition is independent of which orthonormal basis we use. When
we wish to emphasize which metric $\omega$ the norm and inner product
depend on, we will write $|\cdot|_{\omega}$, and $(\cdot,\cdot)_{\omega}.$

The Hodge-star in our setting is defined by the relation \begin{equation}
\alpha\wedge*J(\beta)=(\alpha,\beta)\omega_{n}.\label{eq:hodgestardef}\end{equation}
For an example, if we choose orthonormal coordinates at a point, then
in terms of these we have that \[
*dx_{I}\wedge d\xi_{J}=c_{IJ}\cdot dx_{J^{c}}\wedge d\xi_{I^{c}},\]
for a constant $c_{IJ}=\pm1$ chosen so that \eqref{eq:hodgestardef}
is true; here $I^{c}$ denotes the complementary index of $I$. We
will later investigate the constant $c_{IJ}$ more carefully. 

The integral of an $(n,n)-$form $\alpha=\alpha_{0}(x)c_{n}dx\wedge d\xi$,
is defined by\begin{equation}
\int_{\mathbb{R}^{n}\times\mathbb{R}^{n}}\alpha_{0}(x)c_{n}dx\wedge d\xi=\int_{\mathbb{R}^{n}}\alpha_{0}(x)dx,\label{eq:integralnnform}\end{equation}
 and this gives us an $L^{2}$-structure on the space of forms: \[
L_{p,q}^{2}=\{(p,q)-\text{forms}\,\alpha:\int_{\mathbb{R}^{n}\times\mathbb{R}^{n}}|\alpha|^{2}\omega_{n}<+\infty\}.\]
We will later consider a weighted version of this $L^{2}-$space.
We remark that in defining the integral \eqref{eq:integralnnform}
we have fixed a volume element $d\xi$ on which the integral thus
depends.

\section{Comparison with the complex theory}

In this section, we will consider how super forms correspond to complex
forms. Let us begin in the linear setting, that is, we consider only
forms at a single point, say $x_{0}\in\mathbb{R}^{n}$. Let $\omega$
be an $\mathbb{R}-$K{\"a}hler form. At the point $x_{0}$, we choose
coordinates $(x_{1},...,x_{n},\xi_{1},...,\xi_{n})$ for $\mathbb{R}^{n}\times\mathbb{R}^{n}$
such that \[
\omega(x_{0})=\sum_{k=1}^{n}dx_{k}\wedge d\xi_{k}.\]
Since we will consider complex forms as well, we let $(z_{1},...,z_{n})$
be the standard complex coordinates of $\mathbb{C}^{n}$. We will
use the notation, \[
dV_{i}=dx_{i}\wedge d\xi_{i},\, dV_{i}^{\mathbb{C}}=dz_{i}\wedge d\bar{z}_{i}\]
and for a multi-index $I=(i_{1},...,i_{p})$, we let\[
dV_{I}=dV_{i_{1}}\wedge...\wedge dV_{i_{p}},\, dV_{I}^{\mathbb{C}}=dV_{i_{1}}^{\mathbb{C}}\wedge...\wedge dV_{i_{p}}^{\mathbb{C}}.\]
Now, define \[
\Theta_{I,J,K}=dx_{J}\wedge d\xi_{K}\wedge dV_{I},\]
 for disjoint indices $I,J,$ and $K$. We also define the complex
form \[
\Theta_{I,J,K}^{\mathbb{C}}=dz_{J}\wedge d\bar{z}_{K}\wedge dV_{I}^{\mathbb{C}}.\]
Every super form $\alpha$ can at a fixed point be written as a linear
combination\[
\alpha=\sum\alpha_{I,J,K}\Theta_{I,J,K},\]
where the coefficients $\alpha_{I,J,K}$ are real numbers; we define
a map $\mathcal{C}$ which takes super forms to complex forms by,
\begin{equation}
\mathcal{C}(\alpha)=\sum\alpha_{I,J,K}\Theta_{I,J,K}^{\mathbb{C}}.\label{eq:complexsuperform}\end{equation}
The map $\mathcal{C}$ is linear by definition, and it is also injective,
since $\alpha=0$ is equivalent to $\mathcal{C}(\alpha)=0$. However,
only complex forms of the type \eqref{eq:complexsuperform} with real
coefficients correspond to a super form, and thus, the correspondence
describes an isomorphism between the vector space of super forms at
a fixed point and the vector space of complex forms of the form \eqref{eq:complexsuperform}
with real coefficients $\alpha_{I,J,K}$. This latter space is, from
the complex point of view, not very natural and depends very much
on the choice of coordinates. For instance, a generic change of coordinates
on the complex side does not leave this space invariant. 

The operation of multiplying with the $\mathbb{R}-$K{\"a}hler form $\omega$
is sufficiently important to deserve its own notation: 
\begin{defn}
We define the operator \[
L:\{k-\text{forms}\}\rightarrow\{(k+2)-\text{forms}\}\]
by letting \[
L(\alpha)=\omega\wedge\alpha,\]
where $\alpha$ is a $k-$form. The dual $\Lambda$ of the operator
$L$ is defined by, \[
(L(\alpha),\beta)=(\alpha,\Lambda(\beta)),\]
for $\beta$ a $(k+2)-$form.
\end{defn}
On the complex side, we set the K{\"a}hler form to be $\Omega=\frac{i}{2}\sum_{k=1}^{n}dz_{k}\wedge d\bar{z}_{k}$,
with corresponding operator $L_{\Omega}$. This form $\Omega$ induces
an inner product on the space of complex forms such that the square
of the norm of $\sum\alpha_{I,J,K}\Theta_{I,J,K}^{\mathbb{C}}$ is
equal to $\sum|\alpha_{I,J,K}|^{2}$ in exactly the same way as in
formula \eqref{eq:definitionofnormofform}. The definitions are made
so that, if $\alpha$ is a super form, then the norm of $\mathcal{C}(\alpha)$
measured with respect to $\Omega$, is equal to the norm of $\alpha$
measured with respect to $\omega$. Thus the correspondence $\alpha\leftrightarrow\alpha^{\mathbb{C}}$
is in fact an isometry. We denote by $\Lambda_{\Omega}^{\mathbb{}}$
the dual of $L_{\Omega}$ with respect to to the metric given by $\Omega$.
We have the following:
\begin{prop}
\label{pro:supercomplexcorrespondence}Let $\alpha$ be a $k-$form.
Then \begin{equation}
\mathcal{C}(L\alpha)=\frac{2}{i}L_{\Omega}(\mathcal{C}(\alpha)),\label{eq:complexL}\end{equation}
 \begin{equation}
\mathcal{C}(\Lambda\alpha)=\frac{i}{2}\Lambda_{\Omega}(\mathcal{C}(\alpha)).\label{eq:complexLambda}\end{equation}
 Moreover, \[
\Lambda\alpha=0\Longleftrightarrow\Lambda_{\Omega}\mathcal{C}(\alpha)=0.\]
\end{prop}
\begin{proof}
We let $I+i$ be the multi index $I\cup\{i\}$ and $I-i=I\setminus\{i\}$
which we define to be the empty set if $i\notin I$. First, the formula
\eqref{eq:complexL} is immediate. Next, we claim that \begin{equation}
\Lambda\Theta_{L,M,N}=\sum_{j\in L}\Theta_{L-j,M,N},\label{eq:Lambdappliedtopuretypeform}\end{equation}
where we use the convention that if an index $I,J,$ or $K$ is the
empty set, then $\Theta_{I,J,K}=0$. One realizes this as follows:
we have that \[
L\Theta_{I,J,K}=\sum_{\{i\notin I\bigcup J\bigcup K\}}\Theta_{I+i,J,K}.\]
and that $\Lambda$ is defined by the relation \[
(L(\Theta_{I,J,K}),\Theta_{L,M,N})=(\Theta_{I,J,K},\Lambda(\Theta_{L,M,N})).\]
Here the left hand side is non-zero if and only if there is an $i\notin I\bigcup J\bigcup K,$
such that $I+i=L$ and $J=M$, $K=N$. In this case, the left hand
side is equal to 1, which proves the formula. On the other hand, we
have the well known formula (c.f \cite{Weil}, p.21) \[
\Lambda_{\Omega}\Theta_{L,M,N}^{\mathbb{C}}=\frac{2}{i}\sum_{j\in L}\Theta_{L-j,M,N}^{\mathbb{C}}.\]
 Thus, using linearity of $\Lambda,$ we conclude that formula \eqref{eq:complexLambda}
holds. The last part follows since $\Lambda\alpha=0\Longleftrightarrow\mathcal{C}(\Lambda\alpha)=0\Longleftrightarrow\Lambda_{\Omega}\mathcal{C}(\alpha)=0$.
\end{proof}
The Hodge-star $*_{\Omega}$, acting on complex forms, is defined
by the formula \[
v\wedge*_{\Omega}(\bar{v})=|v|^{2}\Omega_{n}.\]
Let\[
N=\{1,2,...,n\},\]
and recall that we defined \[
c_{p}=(-1)^{p(p-1)/2},\]
for each integer $p$; the number $c_{p}$ was chosen so that \[
dx_{I}\wedge d\xi_{I}=c_{p}\, dx_{i_{1}}\wedge d\xi_{i_{1}}\wedge....\wedge dx_{i_{p}}\wedge d\xi_{i_{p}}\]
where $I=(i_{1},...,i_{p})$. Now, it is well known that (c.f. \cite{Weil},
p.20) that \[
*_{\Omega}dz_{A}\wedge d\bar{z}_{B}\wedge dV_{M}^{\mathbb{C}}=\left[i^{p-q}(-1)^{k(k-1)/2+m}(-2i)^{k-n}\right]dz_{A}\wedge d\bar{z}_{B}\wedge dV_{M^{'}}^{\mathbb{C}},\]
with $M^{'}=N\setminus(A\cup B\cup M)$. However, a small calculations
reveals that\[
*dx\wedge d\xi_{B}\wedge dV_{M}=c_{p}c_{q}(-1)^{p+m+pq}dx_{A}\wedge d\xi_{B}\wedge dV_{M^{'}}.\]
Thus, the real and the complex Hodge stars are related by \begin{equation}
\mathcal{C}(*dx\wedge d\xi_{B}\wedge dV_{M})=i^{n}2^{n-k}(-1)^{n}\cdot(*_{\Omega}dz_{A}\wedge d\bar{z}_{B}\wedge dV_{M}^{\mathbb{C}})\label{eq:differentstars}\end{equation}
since a straightforward calculation shows that \[
\frac{c_{p}c_{q}(-1)^{p+m+pq}}{i^{p-q}(-1)^{k(k-1)/2+m}(-2i)^{k-n}}=i^{n}2^{n-k}(-1)^{n}.\]
Thus, if $\alpha$ is a $k$-form, then \[
\mathcal{C}(*\alpha)=(i^{n}2^{n-k}(-1)^{n})*_{\Omega}(\mathcal{C}(\alpha)).\]
From the complex theory, if $v$ is a complex form, there is a relation
between $*^{\Omega}L_{\Omega}^{r}v$ and $L_{\Omega}^{n-r-k}v$ given
by the following theorem (cf. \cite{Weil}, Theorem 2):
\begin{thm}
If $v=\sum_{|I|=p,|J|=q,|M|=m}v_{I,J,M}\Theta_{I,J,M}^{\mathbb{C}}$,
then 

\[
*_{\Omega}L_{\Omega}^{r}v=i^{p-q}(-1)^{k(k+1)/2}\frac{r!}{(n-k-r)!}L_{\Omega}^{n-r-k}v.\]

\end{thm}
If we apply the above theorem to $\mathcal{C}(\alpha)$ for a $k$-form
\[
\alpha=\sum_{|I|=p,|J|=q,|M|=m}\alpha_{I,J,M}\Theta_{I,J,M},\]
 using \eqref{eq:differentstars} and that $L_{\Omega}^{N}\mathcal{C}(\alpha)=(\frac{i}{2})^{N}\mathcal{C}(L\alpha),$
we obtain,\[
\frac{1}{i^{n}2^{n-k-2r}(-1)^{n}}*((\frac{i}{2})^{r}L{}^{r}\alpha)=i^{p-q}(-1)^{k(k+1)/2}\frac{r!}{(n-k-r)!}(\frac{i}{2})^{n-k-r}L{}^{n-r-k}\alpha,\]
 which gives us \[
*L{}^{r}\alpha=(-1)^{k(k+1)/2+r+q+m}\frac{r!}{(n-k-r)!}L{}^{n-r-k}\alpha.\]
Thus, we have proved the following:
\begin{thm}
\label{thm:starrelatedtoL}If \[
\alpha=\sum_{|I|=p,|J|=q,|M|=m}\alpha_{I,J,M}\Theta_{I,J,M},\]
and $k=p+q+2m$, then 
\end{thm}
\begin{equation}
*L^{r}\alpha=(-1)^{k(k+1)/2+r+q+m}\frac{r!}{(n-k-r)!}L^{n-r-k}\alpha.\label{eq:actualformularelatingstarLrtorL}\end{equation}

This far, we have only compared super forms with complex forms in
the linear setting, that is, at a fixed point. Let us now extend the
map $\mathcal{C}$ to be defined on super forms on all of $\mathbb{R}^{n}$.
Until this point, there has been no need for a relationship between
our real coordinates $(x_{1},.,,,x_{n},\xi_{1},...,\xi_{n})$ and
$(z_{1},.,,,.z_{n})$, but now we make the usual identification $z_{k}=x_{k}+iy_{k}$
for each $k=1,...,n$. For \[
\alpha(x)=\sum\alpha_{IJM}(x)\Theta_{IJM},\]
where $\alpha_{IJM}(\cdot)$ are functions on $\mathbb{R}^{n}$, we
define\[
\mathcal{C}(\alpha)(z)=\sum\alpha_{IJM}(x)\Theta_{IJM}^{\mathbb{C}},\]
where $x=(z+\bar{z})/2.$ 
\begin{prop}
\label{pro:complexvsrealdifferntials}For any super form $\alpha$
we have that \[
\mathcal{C}(d\alpha)=2\partial\mathcal{C}(\alpha),\]
\[
\mathcal{C}(d^{\#}\alpha)=2\bar{\partial}\mathcal{C}(\alpha).\]
\end{prop}
\begin{proof}
Let $\alpha(x)=\sum_{I,J,M}\alpha_{IJM}(x)dx_{I}\wedge d\xi_{J}\wedge dV_{M}.$
Then\[
\mathcal{C}(d\alpha)=\mathcal{C}(\sum_{I,J,M,l}\frac{\partial\alpha_{IJM}(x)}{\partial x_{l}}dx_{l}\wedge dx_{I}\wedge d\xi_{J}\wedge dV_{M})=\]
\[
=\sum_{I,J,M,l}\frac{\partial\alpha_{IJM}(x)}{\partial x_{l}}dz_{l}\wedge dz_{I}\wedge d\bar{z}_{J}\wedge dV_{M}^{\mathbb{C}}.\]
Since $\frac{\partial}{\partial z_{l}}=\frac{1}{2}(\frac{\partial}{\partial x_{l}}-i\frac{\partial}{\partial y_{l}})$,
we see that \[
\frac{\partial\alpha(x)}{\partial x_{l}}=2\frac{\partial\alpha(x)}{\partial z_{l}},\]
and thus\[
\mathcal{C}(d\alpha)=2\partial(\mathcal{C}(\alpha)).\]
The formula for $\bar{\partial}$ follows in the same way. 
\end{proof}
An important formula in complex analysis is the following (c.f \cite{Weil},
p. 42-44):
\begin{thm}
For any complex form $v$,\[
[\Lambda_{\Omega},\partial]v=-i*_{\Omega}\partial*_{\Omega}v.\]

\end{thm}
Let $\alpha$ be a $k-$form. Then, applying the above theorem, we
get\[
[\Lambda_{\Omega},\partial]\mathcal{C}(\alpha)=-i*_{\Omega}\partial*_{\Omega}\mathcal{C}(\alpha).\]
However, by Propositions \ref{pro:supercomplexcorrespondence} and
\ref{pro:complexvsrealdifferntials}, we notice that \[
[\Lambda_{\Omega},\partial]\mathcal{C}(\alpha)=-i\cdot\mathcal{C}([\Lambda,d]\alpha),\]
and, by repeated use of \eqref{eq:differentstars}, keeping in mind
that $d*\alpha$ is a $(2n-k+1)$-form, \[
*_{\Omega}\partial*_{\Omega}\mathcal{C}(\alpha)=\left((i^{n}2^{n-k}(-1)^{n})^{-1}\right)*_{\Omega}\partial\mathcal{C}(*\alpha)=\]
\[
=\left((i^{n}2^{n-k}(-1)^{n}i^{n}2^{n-(2n-k+1)}(-1)^{n})^{-1}\right)\mathcal{C}(*\frac{1}{2}d(*\alpha))=(-1)^{n}\mathcal{C}(*d*\alpha).\]
This gives us that \[
-i(-1)^{n}\mathcal{C}(*d*\alpha)=-i\cdot\mathcal{C}([\Lambda,d]\alpha),\]
and so we arrive at:
\begin{thm}
\label{thm:commutatorrelationfromthecomplexcase}For any form $\alpha$
we have \[
[\Lambda,d]\alpha=(-1)^{n}*d*(\alpha).\]

\end{thm}
Let us conclude this section with some elementary observations:
\begin{lem}
\label{lem:Formlemma}For any $k-$form $\alpha$ we have \begin{equation}
**\alpha=(-1)^{n-k}\alpha\label{eq:starappliedtwice}\end{equation}
\[
\Lambda J\alpha=-J\Lambda\alpha,\]
and \[
*J=(-1)^{n}J*.\]
\end{lem}
\begin{proof}
Since every form $\alpha$ is a linear combination of forms of the
type $dx_{A}\wedge d\xi_{B}\wedge dV_{M}$, we need only to prove
the lemma with $\alpha=dx_{A}\wedge d\xi_{B}\wedge dV_{M}$. One easily
calculates\[
*dx_{A}\wedge d\xi_{B}\wedge dV_{M}=c_{p}c_{q}(-1)^{p+m+qp}dx_{A}\wedge d\xi_{B}\wedge dV_{M'},\]
with $M'=\{1,2,...,n\}\setminus A\cup B\cup M,$ and $p=|A|$, $q=|B|,$$m=|M|$
using the same notation as before. Thus by applying the Hodge-star
twice, the form $dx_{A}\wedge d\xi_{B}\wedge dV_{M}$ will be multiplied
by the constant \[
c_{p}^{2}c_{q}^{2}(-1)^{p+m+qp+p+m'+pq}=(-1)^{m+n-m-p-q}=(-1)^{n-k},\]
where $m'=|M'|=n-p-q-m$, which proves the first formula. The second
formula follows by using \eqref{eq:Lambdappliedtopuretypeform}. The
last formula follows from direct calculations: \[
J*dx_{A}\wedge d\xi_{B}\wedge dV_{M}=c_{p}c_{q}(-1)^{p+m+m'}dx_{A}\wedge d\xi_{B}\wedge dV_{M'},\]

\[
*J(dx_{A}\wedge d\xi_{B}\wedge dV_{M})=c_{p}c_{q}(-1)^{q}dx_{A}\wedge d\xi_{B}\wedge dV_{M'}.\]
Thus $*J$ differs from $J*$ by the constant \[
c_{p}c_{q}(-1)^{p+m+m'}\cdot c_{p}c_{q}(-1)^{q}=(-1)^{n}.\]

\end{proof}
Finally, we note the following corollary of Theorem \ref{thm:commutatorrelationfromthecomplexcase}:
\begin{cor}
\label{cor:duallambdad}For a form $\alpha$, we have, \[
*d^{\#}*\alpha=(-1)^{n+1}[\Lambda,d^{\#}]\alpha.\]
\end{cor}
\begin{proof}
Let us consider the expression $J[\Lambda,d]J\alpha$. From the definition
of $d^{\#}$ we get:\[
J[\Lambda,d]J\alpha=J\Lambda dJ\alpha-Jd\Lambda J\alpha=J\Lambda Jd^{\#}\alpha+JdJ\Lambda\alpha=\]
\[
=-J^{2}\Lambda d^{\#}\alpha+J^{2}d^{\#}\Lambda\alpha=-[\Lambda,d^{\#}]\alpha,\]
 where we used that \[
J\Lambda=-\Lambda J.\]
 On the other hand, Theorem \ref{thm:commutatorrelationfromthecomplexcase}
says that \[
J[\Lambda,d]J\alpha=(-1)^{n}J(*d*)J\alpha,\]
and since $*d*J=J*d^{\#}*$ by applying Lemma \ref{lem:Formlemma},
we have proved that indeed, \[
*d^{\#}*\alpha=(-1)^{n+1}[\Lambda,d^{\#}]\alpha.\]

\end{proof}

\section{Primitive super forms}

In this section we take the opportunity to introduce the notion of
primitivity for super forms and establish expected results by once
again comparing with the complex setting. 
\begin{prop}
\label{pro:commutatorLlambda}Let $\alpha$ be a $(p,q)-$form with
$p+q=k.$ Then \begin{equation}
[\Lambda,L^{s}]\alpha=C_{k,s}L^{s-1}\alpha,\label{eq:commutatoralmbdaLformula}\end{equation}
 with \[
C_{k,s}=s(n-k+1-s).\]
\end{prop}
\begin{proof}
The result follows from the complex theory (c.f \cite{Weil}): if
$v$ is any complex form then \[
[\Lambda_{\Omega},L_{\Omega}^{s}]v=C_{k,s}L_{\Omega}^{s-1}v,\]
and the result follows by letting $v=\mathcal{C}(\alpha)$ and by
repeatedly applying Proposition \ref{pro:supercomplexcorrespondence}. 
\end{proof}
Let us define an important concept in this setting:
\begin{defn}
A form $\alpha$ is called \emph{primitive} if $\alpha$ satisfies
\[
\Lambda\alpha=0.\]

\end{defn}
Note that, in view of Proposition \ref{pro:supercomplexcorrespondence},
$\alpha$ is primitive if and only if $\mathcal{C}(\alpha)$ is primitive
(a complex form $v$ is primitive if $\Lambda_{\Omega}v=0$). The
importance of primitive forms is that they are easier to work with
than just any arbitrary form, combined with the fact that any form
can be decomposed into primitive components in the following sense: 
\begin{prop}
\label{pro:primitivedecompoistion}Let $\alpha$ be a $k$-form. Then
we can write $\alpha$ as \[
\alpha=\alpha_{0}+L\alpha_{1}+...+L^{s}\alpha_{s},\]
 where each $\alpha_{j}$ is a primitive $(k-2j)$-form. Moreover,
the terms of the sum are pairwise orthogonal.\end{prop}
\begin{proof}
The result is well known in the complex case (c.f \cite{Weil}). Thus
we know that the formula holds for $\mathcal{C}(\alpha)$, that is
\[
\mathcal{C}(\alpha)=\alpha_{0}^{'}+L\alpha_{1}^{'}+...+L^{s}\alpha_{s}^{'}\]
 where each $\alpha_{j}^{'}$ is a primitive, complex, $(k-2j)-$form.
But since $\mathcal{C}(\alpha)$ has real coefficients we can assume
that each $\alpha_{j}^{'}$ have real coefficients as well. Thus,
each $\alpha_{j}^{'}$ is in fact $\mathcal{C}(\alpha_{j})$ for some
super form $\alpha_{j}$. By Proposition \ref{pro:supercomplexcorrespondence},
each $\alpha_{j}$ is primitive, and the property of being pairwise
orthogonal is immediate since the correspondence $\alpha\longleftrightarrow\mathcal{C}(\alpha)$
is an isometry. \end{proof}
\begin{prop}
L\label{pro:primitivity}et $\alpha$ be a $k-$form. If \[
L^{n-k+s}\alpha=0\]
then \[
\alpha=\sum_{0\leq j\leq s-1}L^{j}\alpha_{j}\]
with $\alpha_{j}$ a primitive $(k-2j)-$forms. Moreover, $\alpha$
is primitive iff \[
L^{n-k+1}\alpha=0.\]
 \end{prop}
\begin{proof}
The formulas are well known in the complex case and translates into
our setting in the same way as above. 
\end{proof}
The main theorem of this section, known in the complex case as the
Lefschetz isomorphism Theorem, is given by the following:
\begin{thm}
Let $k\leq n$. Then the operator \[
L^{n-k}:\{\mbox{k}-forms\}\rightarrow\{\mbox{(2n-k)}-forms\}\]
is an isomorphism. \end{thm}
\begin{proof}
From the complex setting, we know that $(L_{\Omega})^{n-k}$ is an
isomorphism, and it is easily verified that $k$-forms that are real
linear combinations of $\Theta_{I,J,K}^{\mathbb{C}}$ correspond,
via $L_{\Omega}^{n-k}$, to $(2n-k)$-forms that are real linear combinations
of the same type of degree $(2n-k)$ which establishes the Theorem. 
\end{proof}

\section{$L^{2}$-estimates for the $d$-operato\label{sec:-estimates-for-the}r}

Let us fix a closed, strictly positive, smooth $(1,1)-$form $\omega$.
As we have seen, $\omega$ induces an inner product on the space of
forms, so that for $(p,q)-$forms $\alpha$ and $\beta$, the function
$x\mapsto(\alpha,\beta)(x)$, is a function on $\mathbb{R}^{n}$.
We now define the associated $L^{2}$ inner product: 
\begin{defn}
For $(p,q)-$forms $\alpha$ and $\beta$ we define \[
\bigl\langle\alpha,\beta\bigr\rangle=\int_{\mathbb{R}^{n}\times\mathbb{R}^{n}}(\alpha,\beta)\omega_{n}\]
and we define the associated norm \[
||\alpha||^{2}=\bigl\langle\alpha,\alpha\bigr\rangle.\]
 Observe that, by \eqref{eq:hodgestardef}, we have that\begin{equation}
\bigl\langle\alpha,\beta\bigr\rangle=\int_{\mathbb{R}^{n}\times\mathbb{R}^{n}}(\alpha,\beta)\omega_{n}=\int_{\mathbb{R}^{n}\times\mathbb{R}^{n}}\alpha\wedge*J(\beta).\label{eq:innerpdroductrelationwithhodgestar}\end{equation}
We defined before the space $L_{p,q}^{2}$ as the set of all $(p,q)-$forms
$\alpha$, whose coefficients are integrable, and which satisfies
$||\alpha||^{2}<\infty.$ Moreover, we let \[
L^{2}=\oplus_{p,q=0}^{n}(L_{p,q}^{2}).\]
We consider the operator $d:L_{p,q}^{2}\rightarrow L_{p+1,q}^{2}$
as a closed, densely defined operator, with \[
dom(d)=\{\alpha\in L_{p,q}^{2}:d\alpha\in L_{p+1,q}^{2}\},\]
where $d$ is taken in sense of super currents if $\alpha$ is not
smooth (cf. \cite{Lagerberg}). If $\alpha$ is an $(n,q)$-form it
is understood that $d\alpha=0$. By standard arguments, smooth $(p,q)-$forms
with compact support is dense in $L_{p,q}^{2}$ and each such form
is in $dom(d)$. Thus, $dom(d)$ is indeed dense in $L_{p,q}^{2}$.
We define the dual of the operator $d$ with respect to the inner
product by, the relation \[
\bigl\langle d^{\star}\alpha,\beta\bigr\rangle=\bigl\langle\alpha,d\beta\bigr\rangle\]
for smooth forms $\alpha,\beta\in L^{2}.$ The dual of $d^{\#}$ is
defined analogously.\end{defn}
\begin{prop}
\label{pro:dualofdrelatedtocommutator}For a smooth, compactly supported
form $\alpha$ we have \[
d^{\star}\alpha=[\Lambda,d^{\#}]\alpha,\]
 and \[
(d^{\#})^{\star}\alpha=-[\Lambda,d]\alpha.\]
\end{prop}
\begin{proof}
By Theorem \ref{thm:commutatorrelationfromthecomplexcase} and Corollary
\ref{cor:duallambdad} it is enough to prove that \[
d^{\star}=(-1)^{n+1}*d^{\#}*\]
and \[
(d^{\#})^{\star}=(-1)^{n}*d*.\]
To this end, let $\alpha$ be a $k$-form and $\beta$ a $(k+1)$-form,
both smooth and compactly supported. By \eqref{eq:innerpdroductrelationwithhodgestar}
and Stokes' formula, \[
\bigl\langle d\alpha,\beta\bigr\rangle=\int d\alpha\wedge*J(\beta)=(-1)^{k+1}\int\alpha\wedge d*J(\beta).\]
Since $d*J(\beta)$ is a $(2n-k)-$form, we know from Lemma \ref{lem:Formlemma}
that \[
d*J(\beta)=(-1)^{n-(2n-k)}**d*J(\beta)=(-1)^{n-(2n-k)}*J(*d^{\#}*\beta),\]
 since $dJ=Jd^{\#}$. Thus, \[
\bigl\langle d\alpha,\beta\bigr\rangle=(-1)^{k+1+n-(2n-k)}\int\alpha\wedge**d*J(\beta)=(-1)^{n+1}\bigl\langle\alpha,*d^{\#}*\beta\bigr\rangle,\]
which proves that \[
d^{\star}\beta=(-1)^{n+1}*d^{\#}*\beta.\]
By Corollary \ref{cor:duallambdad}, we see that indeed\[
d^{\star}=[\Lambda,d^{\#}].\]
The second formula of the Proposition follows in the same way, using
Theorem \ref{thm:commutatorrelationfromthecomplexcase}. 
\end{proof}
There are two natural Laplace operators in our setting:
\begin{defn}
For $\alpha$ any smooth and compactly supported form, we define \[
\Box\alpha=dd^{\star}\alpha+d^{\star}d\alpha\]
and\[
\Box^{\#}\alpha=d^{\#}(d^{\#})^{\star}\alpha+(d^{\#})^{\star}d^{\#}\alpha.\]

\end{defn}
Our previous work can now be applied to show that these operators
are in fact equal:
\begin{prop}
F\label{pro:unweightedlaplaciansareequal}or any smooth and compactly
supported form $\alpha$ we have that\[
\Box\alpha=\Box^{\#}\alpha.\]
\end{prop}
\begin{proof}
By Proposition \ref{pro:dualofdrelatedtocommutator}, we obtain\[
\Box\alpha=(d[\Lambda,d^{\#}]\alpha+[\Lambda,d^{\#}]d\alpha)\]
and \[
\Box^{\#}\alpha=-(d^{\#}[\Lambda,d]\alpha+[\Lambda,d]d^{\#}\alpha).\]
Writing out the terms explicitly one immediately concludes that these
expressions are equal. 
\end{proof}
Let us consider {}``twisted'' versions of these Laplacians: 
\begin{defn}
For $\varphi$ a smooth function, we define \[
d_{\varphi}=e^{\varphi}de^{-\varphi},\]
and \[
d_{\varphi}^{\#}=e^{\varphi}d^{\#}e^{-\varphi}.\]
We define the weighted inner product\[
\bigl\langle\alpha,\beta\bigr\rangle_{\varphi}=\int_{\mathbb{R}^{n}\times\mathbb{R}^{n}}(\alpha,\beta)e^{-\varphi}\omega_{n},\]
and let $L_{\varphi}^{2}=\oplus_{0\leq p,q\leq n}(L_{p,q,\varphi}^{2})$
be the space of forms such that \[
||\alpha||_{\varphi}^{2}:=\bigl\langle\alpha,\alpha\bigr\rangle_{\varphi}<\infty.\]
This is easily seen to be a Hilbert space. We will write $L_{\varphi}^{2}(\omega)$
when we wish to emphasize which $\mathbb{R}-$K{\"a}hler metric $\omega$
we are integrating against in defining $L_{\varphi}^{2}$. The dual
of $d$ and $d^{\#}$ with respect to this inner product will be denoted
$d^{*}$ and $(d^{\#})^{*}$ respectively. We can now introduce the
{}``twisted'' Laplacians:\[
\Box_{\varphi}\alpha=dd^{*}\alpha+d^{*}d\alpha\]
 and \[
\Box_{\varphi}^{\#}\alpha=d_{\varphi}^{\#}(d_{\varphi}^{\#})^{*}\alpha+(d_{\varphi}^{\#})^{*}d_{\varphi}^{\#}\alpha.\]

\end{defn}
Our next task is to relate these Laplacians to each other in the spirit
of Proposition \ref{pro:unweightedlaplaciansareequal}. We begin with
the weighted analogue of Proposition \ref{pro:dualofdrelatedtocommutator}:
\begin{prop}
\label{pro:weighteddualsandcommutators}For any smooth, compactly
supported form $\alpha$, the equations \[
d^{*}\alpha=[\Lambda,d_{\varphi}^{\#}]\alpha,\]
and \[
(d_{\varphi}^{\#})^{*}\alpha=-[\Lambda,d]\alpha,\]
are satisfied. \end{prop}
\begin{proof}
Let $\alpha$ be a $(p,q)$-form and $\beta$ a $(p+1,q)$-form. Then,
we compute\[
\bigl\langle d\alpha,\beta\bigr\rangle_{\varphi}=\int(d\alpha,\beta)e^{-\varphi}=\int(\alpha,d^{\star}(e^{-\varphi}\beta))=\int(\alpha,e^{\varphi}d^{\star}(e^{-\varphi}\beta))e^{-\varphi}.\]
By Proposition \ref{pro:dualofdrelatedtocommutator} we know that
$d^{\star}=[\Lambda,d^{\#}].$ Inserting this into the last integral,
we see that \[
\bigl\langle d\alpha,\beta\bigr\rangle_{\varphi}=\int(\alpha,e^{\varphi}[\Lambda,d^{\#}](e^{-\varphi}\beta))e^{-\varphi}=\int(\alpha,[\Lambda,d_{\varphi}^{\#}]\beta)e^{-\varphi},\]
using that $\Lambda$ commutes with the operation of multiplying with
$e^{\varphi}.$ But this means precisely that \[
d^{*}\beta=[\Lambda,d_{\varphi}^{\#}]\beta,\]
which proves the first formula. The second one follows in the same
way. \end{proof}
\begin{thm}
If $\alpha$ is a smooth and compactly supported form, then \begin{equation}
\Box_{\varphi}\alpha=\Box_{\varphi}^{\#}\alpha+[dd^{\#}\varphi,\Lambda]\alpha.\label{eq:commutatoridentity}\end{equation}
\end{thm}
\begin{proof}
We calculate, using Proposition \ref{pro:weighteddualsandcommutators},\[
\Box_{\varphi}\alpha-\Box_{\varphi}^{\#}\alpha=dd^{*}\alpha+d^{*}d\alpha-(d_{\varphi}^{\#}(d_{\varphi}^{\#})^{*}\alpha+(d_{\varphi}^{\#})^{*}d_{\varphi}^{\#}\alpha)=\]
\[
=d[\Lambda,d_{\varphi}^{\#}]\alpha+[\Lambda,d_{\varphi}^{\#}]d\alpha+d_{\varphi}^{\#}[\Lambda,d]\alpha+[\Lambda,d]d_{\varphi}^{\#}\alpha.\]
Since \[
d_{\varphi}^{\#}\alpha=d^{\#}\alpha-d^{\#}\varphi\wedge\alpha,\]
 we get\[
\Box_{\varphi}\alpha-\Box_{\varphi}^{\#}\alpha=\Box\alpha-\Box^{\#}\alpha-d[\Lambda,d^{\#}\varphi]\alpha-[\Lambda,d^{\#}\varphi]d\alpha-d^{\#}\varphi\wedge[\Lambda,d]\alpha-[\Lambda,d](d^{\#}\varphi\wedge\alpha).\]
Here we identify $d^{\#}\varphi$ with the operator sending $\alpha\mapsto d^{\#}\varphi\wedge\alpha$.
By Proposition \ref{pro:unweightedlaplaciansareequal}, we know that
the un-weighted Laplace operators satisfy $\Box\alpha-\Box^{\#}\alpha=0$.
Thus, \[
\Box_{\varphi}\alpha-\Box_{\varphi}^{\#}\alpha=-(d[\Lambda,d^{\#}\varphi]\alpha+[\Lambda,d^{\#}\varphi]d\alpha+d^{\#}\varphi\wedge[\Lambda,d]\alpha+[\Lambda,d](d^{\#}\varphi\wedge\alpha)).\]
Expanding the commutators, we see that\[
\Box_{\varphi}\alpha-\Box_{\varphi}^{\#}\alpha=-(d(\Lambda(d^{\#}\varphi\wedge\alpha)-d(d^{\#}\varphi\wedge(\Lambda\alpha))+\Lambda(d^{\#}\varphi\wedge d\alpha)-d^{\#}\varphi\wedge(\Lambda d\alpha)+\]
\[
+d^{\#}\varphi\wedge(\Lambda d\alpha)-d^{\#}\varphi(d(\Lambda\alpha))+\Lambda(d(d^{\#}\varphi\wedge\alpha)-d(\Lambda(d^{\#}\varphi\wedge\alpha))).\]
Removing the terms which cancel out, we obtain \[
\Box_{\varphi}\alpha-\Box_{\varphi}^{\#}\alpha=d(d^{\#}\varphi\wedge(\Lambda\alpha))-\Lambda(d^{\#}\varphi\wedge d\alpha)-\Lambda(d(d^{\#}\varphi\wedge\alpha)+d^{\#}\varphi\wedge(d(\Lambda\alpha))=\]
\[
=dd^{\#}\varphi\wedge(\Lambda\alpha)-d^{\#}\varphi\wedge d(\Lambda\alpha)-\Lambda(d^{\#}\varphi\wedge d\alpha)-\Lambda(dd^{\#}\varphi\wedge\alpha)+\Lambda(d^{\#}\varphi\wedge d\alpha)+d^{\#}\varphi\wedge(d(\Lambda\alpha))=\]
\[
=dd^{\#}\varphi\wedge(\Lambda\alpha)-\Lambda(dd^{\#}\varphi\wedge\alpha)=[dd^{\#}\varphi,\Lambda]\alpha.\]
Putting everything together, we conclude that \[
\Box_{\varphi}\alpha=\Box_{\varphi}^{\#}\alpha+[dd^{\#}\varphi,\Lambda]\alpha,\]
as desired. \end{proof}
\begin{example}
Let us consider a concrete example of this identity. Let $n=1$ and
let $f$ be a smooth function with compact support. For the weight
function, we choose $\varphi=x^{2}/2$. Then \[
\Box_{\varphi}f=d^{*}df=d^{*}(f'dx)=-f''+xf'\]
and\[
\Box_{\varphi}^{\#}f=(d_{\varphi}^{\#})^{*}d_{\varphi}^{\#}f=(d_{\varphi}^{\#})^{*}[(f'-xf)d\xi]=-f''+f+xf'.\]
Thus \[
\Box_{\varphi}f=\Box_{\varphi}^{\#}f-f,\]
 and we see that in this case $[dd^{\#}\varphi,\Lambda]=-Id$ as predicted
by formula \eqref{eq:commutatoridentity}. 
\end{example}
Let us take the inner product of identity \eqref{eq:commutatoridentity}
against a form smooth, compactly supported form $\alpha$:\[
\left\langle \Box_{\varphi}\alpha,\alpha\right\rangle _{\varphi}=\left\langle \Box_{\varphi}^{\#}\alpha,\alpha\right\rangle _{\varphi}+\left\langle [\Lambda,dd^{\#}\varphi]\alpha,\alpha\right\rangle _{\varphi}\Longleftrightarrow\]
\[
\left\langle dd^{*}\alpha+d^{*}d\alpha,\alpha\right\rangle _{\varphi}=\left\langle d_{\varphi}^{\#}(d_{\varphi}^{\#})^{*}\alpha+(d_{\varphi}^{\#})^{*}d_{\varphi}^{\#}\alpha,\alpha\right\rangle _{\varphi}+\left\langle [dd^{\#}\varphi,\Lambda]\alpha,\alpha\right\rangle _{\varphi}.\]
By the definition of the adjoint, this expression gives us the following
fundamental identity, which should be compared with the classical
Bochner-Kodaira-Nakano identity of complex analysis (c.f \cite{Nakano}):
\begin{thm}
For every smooth, compactly supported form $\alpha$,\begin{equation}
||d^{*}\alpha||_{\varphi}^{2}+||d\alpha||_{\varphi}^{2}=||(d_{\varphi}^{\#})^{*}\alpha||_{\varphi}^{2}+||d_{\varphi}^{\#}\alpha||_{\varphi}^{2}+\left\langle [dd^{\#}\varphi,\Lambda]\alpha,\alpha\right\rangle _{\varphi}.\label{eq:integratedcommutatoridentity}\end{equation}

\end{thm}
By the following fundamental theorem of functional analysis, such
an equality can be used to prove the existence of solutions of the
$d-$equation (c.f. \cite{Hormander}): 
\begin{thm}
\label{thm:Functionalanalysisthm}Let $E$ and $F$ be two Hilbert
spaces, equipped with norms $||\cdot||_{E}$ and $||\cdot||_{F}$
and let $\mathcal{H}$ be a closed subspace of $F$. Let $\mathcal{L}:E\rightarrow F$
be a closed, densely defined operator such that $dom(\mathcal{L}^{*})$
is dense in $F$, and that $Range(\mathcal{L})\subset\mathcal{H}$
If, for each $\alpha\in dom(\mathcal{L}^{*})\cap\mathcal{H}$, the
inequality\begin{equation}
||\mathcal{L}^{*}\alpha||_{F}^{2}\geq c||\alpha||_{E}^{2}\label{eq:Ladjointineq}\end{equation}
is satisfied for some fixed constant $c>0$, then we can find an element
$\beta\in E$ such that \[
\mathcal{L}\beta=\alpha\]
and \[
||\beta||_{E}^{2}\leq c^{-1}||\alpha||_{F}^{2}.\]

\end{thm}
To apply the above theorem to the operator $d$, we thus need to show
an inequality of the type \eqref{eq:Ladjointineq} with $\mathcal{L}=d$.
Let us begin with: 
\begin{prop}
\label{pro:Approxthm}Let $\alpha\in dom(d^{*})\cap dom(d)$. Then,
if the inequality\begin{equation}
||d\beta||_{\varphi}^{2}+||d^{*}\beta||_{\varphi}^{2}\geq c||\beta||_{\varphi}^{2}\label{eq:betainequality}\end{equation}
holds for all smooth, compactly supported forms $\beta$ with $c>0$,
then \begin{equation}
||d\alpha||_{\varphi}^{2}+||d^{*}\alpha||_{\varphi}^{2}\geq c||\alpha||_{\varphi}^{2}.\label{eq:alphainequality}\end{equation}
\end{prop}
\begin{proof}
We can of course assume that $d^{*}\alpha\in L_{p-1,q}^{2}$, since
otherwise the inequality trivially holds, and the condition $\alpha\in dom(d)$
means precisely that $d\alpha\in L_{p+1,q}^{2}$. First, we show that
if the inequality \eqref{eq:betainequality} holds for $\beta\in dom(d^{*})\cap dom(d)$
with compact support, then the desired inequality \eqref{eq:alphainequality}
holds. Indeed, let $\chi_{R}$ be a smooth bump function which is
1 on the ball defined by $\{x\in\mathbb{R}^{n}:|x|\leq R\}$ and vanishes
outside $\{x\in\mathbb{R}^{n}:|x|\leq2R\}$. Then it is easy to see
that $\chi_{R}\cdot\alpha\in dom(d^{*})$. By assumption, we know
that \[
||d(\chi_{R}\cdot\alpha)||_{\varphi}^{2}+||d^{*}(\chi_{R}\cdot\alpha)||_{\varphi}^{2}\geq c||\chi_{R}\cdot\alpha||_{\varphi}^{2}.\]
But $d(\chi_{R}\alpha)=d\chi_{R}\wedge\alpha+\chi_{R}d\alpha$, so
\[
||d(\chi_{R}\cdot\alpha)||_{\varphi}\leq||d\chi_{R}\wedge\alpha||_{\varphi}+||\chi_{R}d\alpha||_{\varphi}.\]
The first term satisfies\[
||d\chi_{R}\wedge\alpha||_{\varphi}^{2}=\int|d\chi_{R}\wedge\alpha|^{2}e^{-\varphi}\omega_{n}\leq\int|d\chi_{R}|^{2}|\alpha|^{2}e^{-\varphi}\omega_{n},\]
and by the assumptions on $\alpha$, this term tend to zero by the
dominated convergence theorem, since $|d\chi_{R}(x)|\rightarrow0$
pointwise as $R\rightarrow0$. Since $d\alpha$ belongs to $L^{2}$,
the second term tends to $||d\alpha||^{2}$, and thus we see that
\[
\lim_{R\rightarrow\infty}||d(\chi_{R}\cdot\alpha)||_{\varphi}^{2}\leq||d\alpha||_{\varphi}^{2}.\]
For the term $d^{*}(\chi_{R}\cdot\alpha)),$ a straightforward calculation
reveals that \[
d^{*}(\chi_{R}\cdot\alpha))=\pm*d^{\#}\chi_{R}\wedge(*\alpha)\pm\chi_{R}d^{*}\alpha,\]
and consequently, \[
||d^{*}(\chi_{R}\cdot\alpha)||_{\varphi}\leq||d^{\#}\chi_{R}\wedge(*\alpha)||_{\varphi}+||\chi_{R}d^{*}\alpha||_{\varphi}.\]
By the same argument as above this implies that \[
\lim_{R\rightarrow\infty}||d^{*}(\chi_{R}\cdot\alpha)||_{\varphi}^{2}\leq||d^{*}\alpha||_{\varphi}^{2}.\]
Combining these observations, we obtain \[
||d\alpha||_{\varphi}^{2}+||d^{*}\alpha||_{\varphi}^{2}\geq c\lim_{R\rightarrow\infty}||\chi_{R}\cdot\alpha||_{\varphi}^{2}=c||\alpha||_{\varphi}^{2},\]
as desired. The proof will thus be complete if we show that the hypothesis
of the proposition implies that the inequality \eqref{eq:alphainequality}
holds for every $\alpha\in dom(d^{*})\cap dom(d)$ with compact support.
But for such an $\alpha$, if we let $\psi_{\epsilon}$ be an approximation
of the identity, it is not hard to show that\[
\alpha*\psi_{\epsilon}\in dom(d^{*}),\]
and moreover, as $\epsilon\rightarrow0$, \[
||d(\alpha*\psi_{\epsilon})||_{\varphi}^{2}\rightarrow||d\alpha||_{\varphi}^{2}.\]
Since $d^{*}$ is a first order differential operator with smooth
coefficients, the same holds true for $d^{*}$ in view of Friedrich's
lemma (c.f. \cite{Friedrichs} or \cite{Hormander} Lemma 1.2.2 which
applies analogously in our setting), that is,\[
||d^{*}(\alpha*\psi_{\epsilon})||_{\varphi}^{2}\rightarrow||d^{*}\alpha||_{\varphi}^{2},\]
 as $\epsilon\rightarrow0$. Thus, since we know that \eqref{eq:betainequality}
holds with $\beta=\alpha*\psi_{\epsilon}$, we see that the inequality
\eqref{eq:alphainequality} holds, and we are done. 
\end{proof}
Now, let $\varphi$ be a smooth convex function such that $dd^{\#}\varphi\geq\epsilon\omega$
for some fixed $\epsilon>0$. We claim that this implies that \[
\left\langle [dd^{\#}\varphi,\Lambda]\alpha,\alpha\right\rangle \geq\epsilon p||\alpha||^{2},\]
for $\alpha$ a $(p,n)-$form. Indeed, by standard linear algebra,
we can at each fixed point $x_{0}$ find orthogonal coordinates in
which $\omega_{x_{0}}=\sum_{i=1}^{n}dx_{i}\wedge d\xi_{i}$ and $dd^{\#}\varphi_{x_{0}}=\sum_{i=1}^{n}\lambda_{i}dx_{i}\wedge d\xi_{i}$
where $\lambda_{i}$ are ordered in such a way that $\lambda_{1}\leq....\leq\lambda_{n}$.
Since this holds for any point $x_{0}$ we can for each $i$ consider
$\lambda_{i}$ as a function on $\mathbb{R}^{n}$ which will depend
continuously on the point $x_{0}$, and by the assumption on $\varphi$
we will have that $\lambda_{n}(x)\geq...\geq\lambda_{1}(x)\geq\epsilon$
for all $x$. For a $(p,n)-$form \[
\alpha=\sum_{|I|=p}\alpha_{I}dx_{I}\wedge d\xi,\]
a calculation reveals (c.f. \cite{Demailly}, p. 69) that the pointwise
inner product at the point $x_{0}$ satisfies \[
([dd^{\#}\varphi,\Lambda]\alpha,\alpha])_{x_{0}}=\sum_{|I|=p}(\sum_{i\in I}\lambda_{i}(x_{0}))|\alpha_{I}|^{2}(x_{0})\geq(\lambda_{1}(x_{0})+...+\lambda_{p}(x_{0}))|\alpha|^{2}(x_{0}).\]
Since $x_{0}$ was arbitrary  we infer that \[
\left\langle [dd^{\#}\varphi,\Lambda]\alpha,\alpha\right\rangle \geq p\epsilon||\alpha||^{2},\]
and thus, by \eqref{eq:integratedcommutatoridentity}, we obtain the
inequality \begin{equation}
||d^{*}\alpha||_{\varphi}^{2}+||d\alpha||_{\varphi}^{2}\geq p\epsilon||\alpha||_{\varphi}^{2},\label{eq:commutatorforddsharpvarphi}\end{equation}
 for every smooth, compactly supported form $\alpha.$ By Proposition
\ref{pro:Approxthm} the above inequality will then hold for every
$\alpha\in dom(d)\cap dom(d^{*}).$ If moreover $d\alpha=0$, this
implies that \[
||d^{*}\alpha||_{\varphi}^{2}\geq p\epsilon||\alpha||_{\varphi}^{2},\]
 and so we can apply Theorem \ref{thm:Functionalanalysisthm} with
$\mathcal{H}=ker(d)$ (which is a closed subspace) to obtain:
\begin{thm}
Let $\omega$ be an $\mathbb{R}-$K{\"a}hler form and let $\varphi$ be
a smooth function such that $dd^{\#}\varphi\geq\epsilon\omega$ for
some $\epsilon>0$. If $\beta\in L_{p,n,\varphi}^{2}$ satisfies that
$d\beta=0$ and if $p\geq1$, then we can find an $\alpha\in L_{p-1,n,\varphi}^{2}$
such that \[
d\alpha=\beta\]
and \[
||\alpha||_{\varphi}^{2}\leq\frac{1}{p\epsilon}||\beta||_{\varphi}^{2}.\]

\end{thm}
If we instead let the $\mathbb{R}$-K{\"a}hler form $\omega$ be given
by $\omega=dd^{\#}\varphi$ for some smooth, convex function $\varphi$,
then Proposition \ref{pro:commutatorLlambda} tells us that $[dd^{\#}\varphi,\Lambda]\alpha=(k-n)\alpha$,
if $\alpha$ is a $k$-form. Thus, \begin{equation}
||d^{*}\alpha||_{\varphi}^{2}+||d\alpha||_{\varphi}^{2}\geq(k-n)||\alpha||_{\varphi}^{2},\label{eq:commutatorforddsharpvarphi-1}\end{equation}
 for every smooth, compactly supported form $\alpha$ and we get the
following result:
\begin{thm}
Assume that $\omega=dd^{\#}\varphi>0$ for a smooth convex function
$\varphi$. If $\beta\in L_{p,q,\varphi}^{2}$ is a $k-$form with
$k>n$ and $p\geq1$ such that $d\beta=0$, then we can find a $\alpha\in L_{p-1,q,\varphi}^{2}$
such that \[
d\alpha=\beta\]
and \[
||\alpha||_{\varphi}^{2}\leq\frac{1}{k-n}||\beta||_{\varphi}^{2}.\]

\end{thm}
Now, if $\tilde{\beta}$ is a closed $(p,n)$-form, we can by the
virtue of the above theorem solve the equation \[
d\tilde{\alpha}=\tilde{\beta}\]
with the estimate \[
||\tilde{\alpha}||_{\varphi}^{2}\leq\frac{1}{p}||\tilde{\beta}||_{\varphi}^{2}.\]
Let us write this in coordinates: if $\tilde{\alpha}=\alpha\wedge d\xi$
with $\alpha=\sum_{|K|=p-1}\alpha_{K}dx_{K},$ and $\tilde{\beta}=\beta\wedge d\xi$
with $\beta=\sum_{|L|=p}\beta_{L}dx_{L},$ then \[
\int|\alpha|_{dd^{\#}\varphi}^{2}c_{n}dx\wedge d\xi\leq\frac{1}{p}\int|\beta|_{dd^{\#}\varphi}^{2}c_{n}dx\wedge d\xi.\]
 This follows since $|\tilde{\alpha}|_{dd^{\#}\varphi}^{2}=|\alpha|_{dd^{\#}\varphi}^{2}\det(\omega_{ij})^{-1},$
and similarly for $|\tilde{\beta}|_{dd^{\#}\varphi}^{2}$. Let us
reverse this argument: if $\beta$ is a closed $(p,0)$-form such
that $\int|\beta|_{dd^{\#}\varphi}^{2}dx\wedge d\xi<+\infty$, we
can consider the $(p,n)$ form $\tilde{\beta}=\beta\wedge d\xi$,
which will also be closed. Then \[
|\tilde{\beta}|_{dd^{\#}\varphi}^{2}=|\beta|_{dd^{\#}\varphi}^{2}\det(\omega_{ij})^{-1},\]
and by the above we have that $\tilde{\beta}\in L_{\varphi}^{2}$.
Thus we can solve $d\tilde{\alpha}=\tilde{\beta}$, for some $(p-1,n)$-
form $\tilde{\alpha}$. But, $\tilde{\alpha}=\alpha\wedge d\xi$ for
some $(p-1,0)$-form $\alpha,$ and we must have $d\alpha=\beta.$
Thus we arrive at:
\begin{thm}
\label{thm:L2existancetheoremforL2forms}For a closed $(p,0)$-form
$\beta$ such that $\int_{\mathbb{R}^{n}}|\beta|_{dd^{\#}\varphi}^{2}e^{-\varphi}dx<\infty$,
we can solve \[
d\alpha=\beta,\]
with \begin{equation}
\int_{\mathbb{R}^{n}}|\alpha|_{dd^{\#}\varphi}^{2}e^{-\varphi}dx\leq\frac{1}{p}\int_{\mathbb{R}^{n}}|\beta|_{dd^{\#}\varphi}^{2}e^{-\varphi}dx.\label{eq:L2estimatep0forms}\end{equation}

\end{thm}
It is interesting to note that when $p=1$ the left-hand-side of \eqref{eq:L2estimatep0forms}
does not depend on the K{\"a}hler metric $dd^{\#}\varphi$. 
\begin{rem}
Let us explain briefly how our setting is related to solving the $\overline{\partial}-$equation
on a holomorphic line bundle $L$ over a compact K{\"a}hler manifold $X$
(see \cite{Demailly} or \cite{Berndtsson} for a detailed account).
Let $\varphi$ be a metric on $L$, inducing a hermitian structure
on $L$, and let $\nabla$ be the Chern connection of $L$. Strictly
speaking, $\varphi$ is a collection of smooth functions $\{\varphi_{i}\}$,
each defined on an open set of $U_{i}\subset X$ corresponding to
a trivialization of $L$. If $s\in H^{0}(X,L)$, then the norm of
$s$ is locally given by $x\mapsto|s_{i}|^{2}(x)e^{-\varphi_{i}(x)}$,
where $s_{i}$ is a local representative of $s$ using the trivialization
of $L$. We can thus perceive $|s|e^{-\varphi}$ as a globally defined
function on $X$. As is well known, we can write the connection $\nabla$
as $\nabla=\nabla'+\nabla''$, where $\nabla''=\overline{\partial}$,
and we can consider the duals of $\nabla'$ and $\nabla''$ with respect
to the metric $\varphi$, and denote them by $(\nabla')^{*}$ and
$(\nabla'')^{*}$. Then the classical Bochner-Kodaira-Nakano identity
states that \[
\nabla'(\nabla')^{*}+(\nabla')^{*}\nabla'=\nabla''(\nabla'')^{*}+(\nabla'')^{*}\nabla''+[dd^{c}\varphi,\Lambda],\]
where one can show that $dd^{c}\varphi$ is the curvature operator
associated with $\nabla.$ In the same way as in this article, this
identity can be used to show the solvability of the $\bar{\partial}-$equation
(the argument is basically due to H{\"o}rmander \cite{Hormander}): \end{rem}
\begin{thm}
\emph{(H{\"o}rmander}, \cite{Hormander}\emph{)} Let $\beta$ be a $(p,q)-$form
with values on $L$ such that $\overline{\partial}\beta=0$ and $\int_{X}|\beta|^{2}e^{-\varphi}dV_{X}<+\infty$,
and assume we can find a metric $\varphi$ on $L$ such that $([dd^{c}\varphi,\Lambda]\alpha,\alpha)\geq c||\alpha||^{2}$
for each compactly supported $(p,q+1)$-form $\alpha$ with values
in $L$. Then we can solve the equation \[
\overline{\partial}\alpha=\beta\]
with \[
\int_{X}|\alpha|^{2}e^{-\varphi}dV_{X}\leq\frac{1}{c}\int_{X}|\beta|^{2}e^{-\varphi}dV_{X},\]
where $dV_{X}$ is the volume element on $X$. 
\end{thm}
Now, let us instead consider the Laplace operator $\Box_{\varphi}$:
if $k>n$ and $\alpha$ is a $k$-form in $L_{\varphi}^{2}$ we know
that \eqref{eq:commutatorforddsharpvarphi} holds, that is \[
||d\alpha||^{2}+||d^{*}\alpha||^{2}\geq(k-n)||\alpha||^{2},\]
under the assumption that the metric in question is $dd^{\#}\varphi$.
As we already have seen, the left-hand-side is equal to $\left\langle \Box_{\varphi}\alpha,\alpha\right\rangle ,$
and by the Cauchy-inequality\[
\left\langle \Box_{\varphi}\alpha,\alpha\right\rangle \leq||\Box_{\varphi}\alpha||\cdot||\alpha||.\]
Thus\[
||\Box_{\varphi}\alpha||\geq(k-n)||\alpha||.\]
 Since $\Box_{\varphi}$ is self-adjoint we can apply Theorem \ref{thm:Functionalanalysisthm}
to obtain:
\begin{thm}
With the notation and assumptions above, we can for each $d-$closed
$k$-form $\beta$ solve the equation \[
\Box_{\varphi}\alpha=\beta,\]
with \[
||\alpha||_{\varphi}^{2}\leq(k-n)^{2}||\beta||_{\varphi}^{2}.\]

\end{thm}

\section{The Legendre transform}

We recall the definition and some properties of the Legendre transform:
\begin{defn}
Let $f$ be a convex function on $\mathbb{R}^{n}$. The Legendre transform
of $f$ is given by\begin{equation}
f^{*}(y)=\sup_{x\in\mathbb{R}^{n}}(x\cdot y-f(x)),\label{eq:Legendre}\end{equation}
for $y\in\mathbb{R}^{n}$. 
\end{defn}
The Legendre transform of a convex function is again convex, and $(f^{*})^{*}=f$.
Let us assume that $f$ is smooth. Then the supremum in \eqref{eq:Legendre}
(if it is not equal to $+\infty$, which we always shall assume in
this section) is achieved at a point $x(y)$ for which \[
y=\nabla f(x(y)),\]
where $\nabla f$ denotes the gradient of $f$; to see this is simply
a matter of differentiating the expression inside of the supremum.
Thus \[
f^{*}(y)=x(y)\cdot y-f(x(y)).\]
By a small calculation this implies that if $y=\nabla f(x(y))$ as
above, then \[
\nabla f^{*}(y)=x.\]
Thus, if we consider the map $\psi:\mathbb{R}^{n}\rightarrow\mathbb{R}^{n}$,
given by \[
\psi(x)=\nabla f(x),\]
then \[
\psi^{-1}(y)=\nabla f^{*}(y).\]
For any smooth map $\psi:\mathbb{R}^{n}\rightarrow\mathbb{R}^{n}$,
the pullback, $\psi^{\star}$, of a form $(p,0)$-form is just the
regular pullback of a differential $p$-form, and we extend $\psi^{\star}$
to act on $(p,q)$-forms by requiring it to be $J$-linear, that is\[
\psi^{\star}J(\alpha)=J(\psi^{\star}\alpha),\]
for any $(p,q)-$form $\alpha.$ This makes sense since, in coordinates,
this makes for \[
\psi^{\star}(\alpha_{IJ}(x)dx_{I}\wedge d\xi_{J})=\alpha_{IJ}(\psi(x))\psi^{\star}dx_{I}\wedge J(\psi^{\star}dx_{J}),\]
 and every $(p,q)$-form can be written as a linear combination of
such forms. 
\begin{prop}
\label{pro:changeofvariablesinintegrals}Let $\psi:\mathbb{R}^{n}\rightarrow\mathbb{R}^{n}$
be a diffeomorphism. Then any integrable $(n,n)-$form $\alpha$ satisfies,
\[
\int_{\mathbb{R}^{n}\times\mathbb{R}^{n}}\alpha=\int_{\mathbb{R}^{n}\times\mathbb{R}^{n}}\psi^{\star}(\alpha)/det(D\psi).\]
\end{prop}
\begin{proof}
This is a simple consequence of the usual change of variable formula
for $n$-forms on $\mathbb{R}^{n}$: let $\alpha=\alpha_{0}(x)c_{n}dx\wedge d\xi$.
Then $\psi^{\star}\alpha(x)=(\alpha_{0}\circ\psi)(x)\cdot(det(D\psi)(x))^{2}c_{n}dx\wedge d\xi.$
Thus \[
\int_{\mathbb{R}^{n}\times\mathbb{R}^{n}}\psi^{\star}(\alpha)/det(D\psi)=\int_{\mathbb{R}^{n}\times\mathbb{R}^{n}}(\alpha_{0}\circ\psi)(x)\cdot det(D\psi)(x)c_{n}dx\wedge d\xi=\]
\[
=\int_{\mathbb{R}^{n}}(\alpha_{0}\circ\psi)(x)\cdot det(D\psi)(x)dx=\int_{\mathbb{R}^{n}}\alpha_{0}(x)dx=\int_{\mathbb{R}^{n}\times\mathbb{R}^{n}}\alpha.\]

\end{proof}
Now, let $\varphi$ be a smooth, strictly convex function, and associate
to $\varphi$ the $\mathbb{R}-$K{\"a}hler form \[
\omega^{\varphi}=dd^{\#}\varphi=\sum_{i,j=1}^{n}\varphi_{ij}dx_{i}\wedge d\xi_{j},\]
with $\varphi_{ij}=\frac{\partial\varphi}{\partial x_{i}\partial x_{j}}$.
Let also $\psi=\nabla\varphi^{*}$ which is a diffeomorphism. Since
$\nabla\varphi^{*}$ and $\nabla\varphi$ are inverse to each other
by the above, we have $(\varphi_{ij}^{*})=(\varphi^{ij})$ where $(\varphi^{ij})$
denotes the inverse of the matrix $(\varphi_{ij})$. Thus, since $\psi^{\star}dx_{i}=\sum_{i,k=1}^{n}\varphi_{ik}^{*}dx_{k}$
and $\psi^{\star}d\xi_{i}=\sum_{i,k=1}^{n}\varphi_{ik}^{*}d\xi_{k}$,
we conclude that\[
\psi^{\star}\omega^{\varphi}=\sum_{i,j,k,l=1}^{n}\varphi_{ij}\varphi^{ik}\varphi^{jl}dx_{k}\wedge d\xi_{l}=\sum_{j,l=1}^{n}\varphi^{jl}dx_{j}\wedge d\xi_{l}.\]
Here we used that the matrix $(\varphi_{ij})$ is symmetric so that
$\sum_{i=1}^{n}\varphi_{ij}\varphi^{ik}=\delta_{jk}$. Thus, we obtain:
\begin{equation}
\psi^{\star}\omega^{\varphi}=\omega^{\varphi^{*}}.\label{eq:LegendreofKahlerform}\end{equation}
Recall from section \ref{sec:Preliminaries} that the norm of a $(p,0)$-form
$\alpha$ satisfies the relation \[
|\alpha|_{dd^{\#}\varphi}^{2}\omega_{n}^{\varphi}=c_{p}\,\alpha\wedge J(\alpha)\wedge\omega_{n-p}^{\varphi}.\]
Applying $\psi^{\star}$ to both sides of this equality and using
\eqref{eq:LegendreofKahlerform}, give us \[
\psi^{\star}(|\alpha|_{dd^{\#}\varphi}^{2}\omega_{n}^{\varphi})=c_{p}\,\psi^{\star}\alpha\wedge J(\psi^{\star}\alpha)\wedge\omega_{n-p}^{\varphi^{*}}.\]
This in turn is equivalent to \[
|\alpha|_{dd^{\#}\varphi}^{2}(\psi(x))\omega_{n}^{\varphi^{*}}=|\psi^{\star}\alpha|_{dd^{\#}\varphi^{*}}^{2}(x)\omega_{n}^{\varphi^{*}},\]
and we have proved:
\begin{prop}
\label{pro:pullbackofp0forms}Let $\alpha$ be a $(p,0)-$form. Then,
at any point $x$, \[
|\psi^{\star}\alpha|_{dd^{\#}\varphi^{*}}^{2}(x)=|\alpha|_{dd^{\#}\varphi}^{2}(\psi(x)).\]

\end{prop}
Let us consider the integral\[
\int|\alpha|^{2}e^{-\varphi}dx\wedge d\xi;\]
we shall see how this integral transform under the Legendre transform
of $\varphi$ under the additional assumption that $\varphi$ is $r-$homogeneous,
that is, when for each $y\in\mathbb{R}^{n},$ \begin{equation}
\varphi(ty)=t^{r}\varphi(y),\label{eq:homogneity}\end{equation}
for $t\geq0$. Differentiating the relation \eqref{eq:homogneity}
with respect to $t$, and evaluating at $t=1$ tells us that, \[
y(\nabla\varphi)(y)=r\varphi(y).\]
This result is sometimes referred to as Euler's theorem on homogeneous
functions. Moreover, we know that \[
\varphi^{*}(x)=y\cdot(\nabla\varphi)(y)-\varphi(y)\]
where $y$ is such that $x=\nabla\varphi(y)$. Thus, \[
\varphi^{*}(x)=(r-1)\varphi(y).\]
Furthermore, \[
\psi^{\star}(dx\wedge d\xi)=det(\varphi_{ij}^{*})^{2}dx\wedge d\xi.\]
Thus, by Proposition \ref{pro:changeofvariablesinintegrals} and Proposition
\ref{pro:pullbackofp0forms}, we obtain \[
\int_{\mathbb{R}^{n}\times\mathbb{R}^{n}}|\alpha|_{\omega^{\varphi}}^{2}e^{-\varphi}c_{n}dx\wedge d\xi=\int_{\mathbb{R}^{n}\times\mathbb{R}^{n}}\psi^{\star}(|\alpha|_{\omega^{\varphi}}^{2}e^{-\varphi}c_{n}dx\wedge d\xi)/det(D\psi)=\]
\[
=\int_{\mathbb{R}^{n}\times\mathbb{R}^{n}}|\psi^{\star}\alpha|_{\omega^{\varphi^{*}}}^{2}e^{-\frac{\varphi^{*}}{r-1}}det(\varphi_{ij}^{*})c_{n}dx\wedge d\xi=\int_{\mathbb{R}^{n}\times\mathbb{R}^{n}}|\psi^{\star}\alpha|_{\omega^{\varphi^{*}}}^{2}e^{-\frac{\varphi^{*}(y)}{r-1}}\omega_{n}^{\varphi^{*}}.\]
Let us record this result as a proposition:
\begin{prop}
\label{pro:Homogeneousintegralcorrespondance}Let $\varphi$ be a
$r$-homogeneous, convex and smooth function, and let $\psi=\nabla\varphi^{*}$.
Then any integrable $(p,0)$-form $\alpha$ satisfies\begin{equation}
\int_{\mathbb{R}^{n}\times\mathbb{R}^{n}}|\alpha|_{\omega^{\varphi}}^{2}e^{-\varphi}c_{n}dx\wedge d\xi=\int_{\mathbb{R}^{n}\times\mathbb{R}^{n}}|\psi^{*}\alpha|_{\omega^{\varphi^{*}}}^{2}e^{-\frac{\varphi^{*}(y)}{r-1}}\omega_{n}^{\varphi^{*}}.\label{eq:integralunderlegendretransform}\end{equation}

\end{prop}
Under these circumstances we can prove the following:
\begin{prop}
\label{pro:solvingdequationwithhomogeneouswieght}Let $\beta\in L_{\phi}^{2}(dd^{\#}\phi)$
be a $(p,0)$-form, with $\phi$ an $r$-homogeneous, convex and smooth
function, with $r>1$. Then we can solve\[
d\alpha=\beta\]
with the estimate\begin{equation}
\int_{\mathbb{R}^{n}\times\mathbb{R}^{n}}|\alpha|_{\omega^{\phi}}^{2}e^{-\phi}\omega_{n}^{\phi}\leq\frac{1}{p(r-1)}\int_{\mathbb{R}^{n}\times\mathbb{R}^{n}}|\beta|_{\omega^{\phi}}^{2}e^{-\phi}\omega_{n}^{\phi}.\label{eq:integralestimateswithlegendretrans}\end{equation}
 \end{prop}
\begin{proof}
If $s=\frac{r}{r-1}$, then it is well known that $\phi^{*}$ is an
$s$-homogeneous function. Thus, if we let \[
\varphi=(r-1)^{1-r}\phi^{*},\]
then $\varphi$ is $s-$homogeneous, and satisfies $\frac{\varphi^{*}}{r-1}=\phi$.
To simplify notation we let $\tau=(r-1)^{1-r}$. If $\alpha$ is a
$(p,0)$-form, then \[
|\alpha|_{\omega^{\varphi}}^{2}=\tau^{-p}|\alpha|_{\omega^{\phi^{*}}}^{2},\]
\[
|\psi^{\star}\alpha|_{\omega^{\varphi^{*}}}^{2}=(r-1)^{-p}|\psi^{\star}\alpha|_{\omega^{\phi}}^{2},\]
\[
\omega_{n}^{\varphi^{*}}=(r-1)^{n}\omega_{n}^{\phi};\]
thus, with $\psi=\nabla\varphi^{*}$ formula \eqref{eq:integralunderlegendretransform}
transforms into \begin{equation}
\int_{\mathbb{R}^{n}\times\mathbb{R}^{n}}|\alpha|_{\omega^{\phi^{*}}}^{2}e^{-\tau\phi^{*}}c_{n}dx\wedge d\xi=(r-1)^{n-pr}\int_{\mathbb{R}^{n}\times\mathbb{R}^{n}}|\psi^{\star}\alpha|_{\omega^{\phi}}^{2}e^{-\phi}\omega_{n}^{\phi}.\label{eq:integralrelationwithlegendretransform}\end{equation}
 Now let us consider the form $\gamma=(\psi^{-1})^{*}\beta$. Then,
since $\psi^{*}(\psi^{-1})^{*}\beta=\beta,$ we see that $\psi^{*}\gamma\in L_{\phi}^{2}(dd^{\#}\phi)$,
and by the above the form $\gamma$ satisfies \[
\int_{\mathbb{R}^{n}\times\mathbb{R}^{n}}|\gamma|_{\omega^{\phi^{*}}}^{2}e^{-\tau\phi^{*}}c_{n}dx\wedge d\xi<+\infty.\]
 Moreover, $\gamma$ is $d$-closed, since $d\psi^{*}=\psi^{*}d$.
Thus, by Theorem \ref{thm:L2existancetheoremforL2forms} we can find
a $(p-1,0)$- form $\eta$ such that\[
d\eta=\gamma\]
and\[
\tau^{-(p-1)}\int_{\mathbb{R}^{n}\times\mathbb{R}^{n}}|\eta|_{\omega^{\phi^{*}}}^{2}e^{-\tau\phi^{*}}c_{n}dx\wedge d\xi\leq\]
\[
\leq\frac{1}{p}\tau^{-p}\int|\gamma|_{\omega^{\phi^{*}}}^{2}e^{-\tau\phi^{*}}c_{n}dx\wedge d\xi,\]
where we used that \[
|\alpha|_{\omega^{\tau\phi^{*}}}^{2}=\tau^{-p}|\alpha|_{\omega^{\phi^{*}}}^{2},\]
for any $(p,0)-$form $\alpha.$ If we apply formula \eqref{eq:integralrelationwithlegendretransform}
to this inequality, with $\alpha=\psi^{\star}\eta$, then since $\psi^{\star}\gamma=\beta$
we see that \[
\int_{\mathbb{R}^{n}\times\mathbb{R}^{n}}|\alpha|_{\omega^{\phi}}^{2}e^{-\phi}\omega_{n}^{\phi}\leq\frac{1}{p(r-1)}\int_{\mathbb{R}^{n}\times\mathbb{R}^{n}}|\beta|_{\omega^{\phi}}^{2}e^{-\phi}\omega_{n}^{\phi}.\]
Also, \[
d\alpha=d\psi^{\star}\eta=\psi^{\star}\gamma=\beta,\]
which concludes the proof. 
\end{proof}
While the estimate of the above proposition is similar to that of
Theorem \ref{thm:L2existancetheoremforL2forms}, it does not seem
to follow from the previous formalism in a direct way. 

From section \ref{sec:-estimates-for-the} we know that in order to
obtain existence theorems for the $d-$operator we need to examine
the commutator term $[dd^{\#}\varphi,\Lambda]$, and we used in that
section the choice $\omega=dd^{\#}\varphi$ for $\varphi$ a smooth,
strictly convex function, to obtain $[dd^{\#}\varphi,\Lambda]\alpha=(k-n)\alpha$,
for $\alpha$ a $k$-form. Let us instead assume that $\mbox{\ensuremath{\phi}\ }$is
a smooth, strictly \emph{concave }function. Then $-dd^{\#}\phi$ is
a closed positive $(1,1)$- form and we can let $\omega^{\phi}=-dd^{\#}\phi$
be our K{\"a}hler form. In this situation we see that \[
[dd^{\#}\phi,\Lambda]\alpha=(n-k)\alpha,\]
 and thus, we have \[
||d^{*}\alpha||_{\phi}^{2}+||d\alpha||_{\phi}^{2}\geq(n-k)||\alpha||_{\phi}^{2},\]
for every smooth, compactly supported form $\alpha$, where $||\alpha||_{\phi}^{2}=\int_{\mathbb{R}^{n}\times\mathbb{R}^{n}}|\alpha|_{\omega^{\phi}}^{2}e^{-\phi}\omega_{n}^{\phi},$
and the dual of $d$ is with respect to this norm. Applying Theorem
\ref{thm:Functionalanalysisthm} we obtain the following:
\begin{prop}
Let $\phi$ be a smooth, strictly concave function. Then, for every
closed $(p,0)$-form $\beta\in L_{\phi}^{2}$, we can solve the equation
\[
d\alpha=\beta,\]
with \[
\int_{\mathbb{R}^{n}\times\mathbb{R}^{n}}|\alpha|_{\omega^{\phi}}^{2}e^{-\phi}\omega_{n}^{\phi}\leq\frac{1}{n-p}\int_{\mathbb{R}^{n}\times\mathbb{R}^{n}}|\beta|_{\omega^{\phi}}^{2}e^{-\phi}\omega_{n}^{\phi}.\]

\end{prop}
It is interesting to compare this result with Proposition \ref{pro:solvingdequationwithhomogeneouswieght}:
let $\varphi$ be $2-$homogeneous and assume that $\beta\in L_{\varphi}^{2}\cap L_{\phi}^{2}.$
Then we can find solutions $\alpha_{1}$ and $\alpha_{2}$ to $d\alpha=\beta$
such that \[
\int_{\mathbb{R}^{n}\times\mathbb{R}^{n}}|\alpha_{1}|_{\omega^{\varphi}}^{2}e^{-\varphi}\omega_{n}^{\varphi}\leq\frac{1}{p}\int_{\mathbb{R}^{n}\times\mathbb{R}^{n}}|\beta|_{\omega^{\varphi}}^{2}e^{-\varphi}\omega_{n}^{\varphi}\]
 \[
\int_{\mathbb{R}^{n}\times\mathbb{R}^{n}}|\alpha_{2}|_{\omega^{\phi}}^{2}e^{-\phi}\omega_{n}^{\phi}\leq\frac{1}{n-p}\int_{\mathbb{R}^{n}\times\mathbb{R}^{n}}|\beta|_{\omega^{\phi}}^{2}e^{-\phi}\omega_{n}^{\phi}.\]
Thus, we can solve the equation $d\alpha=\beta$ with fundamentally
different estimates on the solutions: in one case the weight $\varphi$
is convex and in the other $\phi$ is concave. 

\bibliographystyle{plain}
\nocite{*}
\bibliography{theBochnerNakanoIdentityBib}

\end{document}